\input amstex
\magnification=1200
\documentstyle{amsppt}
\NoRunningHeads
\NoBlackBoxes
\define\Mat{\operatorname{Mat}}
\topmatter
\title Representative dynamics\endtitle
\author Denis V. Juriev\endauthor
\affil ul.Miklukho-Maklaya 20-180, Moscow 117437 Russia\linebreak
(e-mail: denis\@juriev.msk.ru)\endaffil
\date math.OC/9911040\newline November, 07, 1999\enddate
\abstract\nofrills This short note is devoted to the representative dynamics,
which realizes a link between the theory of controlled systems and 
representation theory. Dynamical inverse problem of representation theory for 
controlled systems is considered: to solve it means to correspond a 
representative dynamics to the controlled system.
\endabstract
\endtopmatter
\document
This short note is devoted to the representative dynamics, a new concept, 
which realizes a link between the theory of controlled systems and 
representation theory, and the dynamical inverse problem of representation
theory for controlled systems. Thus, it may be considered as a development
of ideas of the earlier article [1], which follows the general ideology of
inverse problems of representation theory [2].

\subhead 1. Representative dynamics\endsubhead

\definition{Definition 1} Let $\bold X=\bold X(t)=(X_1(t),\ldots X_m(t))$
($X_i(t)\in\Mat_n(\Bbb C)$) be the time-dependent vector of $m$ complex 
$n\times n$ matrices. The {\it representative dynamics\/} is a controlled
system (with constraints) of the form
$$\dot\bold X(t)=F(\bold X(t),a(t))\tag1$$
with the {\sl fixed\/} initial data $\bold X(t_0)$, where the control 
parameter $a(t)=(\frak A(t),\bold e(t))$ is the pair of {\sl any\/} 
associative algebra $\frak A(t)$ from the fixed class of such algebras 
$\Bbb A$, $\bold e(t)=(e_1(t),\ldots e_m(t))$ is {\sl any\/} set of algebraic 
generators of the algebra $\frak A$ (one may claim $\bold e(t)$ to be an 
algebraic basis in $\frak A$) such that the mapping $e_i(t)\mapsto X_i(t)$ 
may be extended to the representation $T(t):\frak A(t)\mapsto\Mat_n(\Bbb C)$ 
of the algebra $\frak A(t)$ in the matrix algebra $\Mat_n(\Bbb C)$ (this is a 
constraint on the control $a(t)$).
\enddefinition

Certainly, the claim that (1) is a representative dynamics restricts the
choice of the function $F$ and initial data $\bold X(t_0)$ because for
each moment $t$ {\sl any\/} admissible choice of the pair $(\frak A(t),
\bold e(t))$ should provide that the set of admissible pairs will not be 
empty in future. 

\remark{Exercise} To describe all representation dynamics with $m\!=\!3$, 
$n\!=\!2$ and the class $\Bbb A$ of all associative algebras $\frak A$ 
with quadratic relations, which are isomorphic as linear spaces to the
symmetric algebra $S^{\cdot}(V)$ over the linear space $V$ spanned by the
fixed elements $e_i$ of $\frak A$ under the Weyl symmetrization mapping.
\endremark

\remark{Remark 1} Let us consider the following equivalence on the set 
$\Cal A$ of all admissible $a=(\frak A,\bold e)$. The pairs $a_1=(\frak A_1,
\bold e_1)$ and $a_2=(\frak A_2,\bold e_2)$ will be equivalent iff the 
algebras $A_1$ and $A_2$ are isomorphic under an isomorphism which maps the 
linear space $V_1$ spanned by the elements of $\bold e_1$ onto the linear 
space $V_2$ spanned by the elements of $\bold e_2$. Then the equivalence
divides the time interval $[t_0,t_1]$, on which the representative dynamics is 
considered, onto the subsets, on which the pairs $a(t)$ are equivalent.
\endremark

\remark{Remark 2} Representative dynamics combines structural and functional
features. The first are accumulated in the class $\Bbb A$ and the least are
expressed by the function $F$. Both are interrelated. The situation is similar
to one in the group theory of special functions [3]. However, the difference 
is essential: in the representative dynamics the functional aspects are not 
derived from the structural ones and hve an independent origin.
\endremark

\subhead 2. Dynamical inverse problem of representation theory for
controlled systems\endsubhead

In the article [1] dynamical inverse problem of representation theory
was considered (see also the review [2] on the general ideology of
inverse problems of representation theory). Below this concept will be
adapted for the controlled systems. The representative dynamics will play
a crucial role.

\definition{Definition 2} Let 
$$\dot x=\varphi(x,u),\tag2$$
be the controlled system, where $x$ is the time-dependent $m$-dimensional
complex vector and $u$ is the control parameter. {\it Dynamical inverse 
problem of representation theory for the controlled system\/} (2) is to 
construct a representative dynamics
$$\dot\bold X=F(\bold X,a)$$
and the function
$$a=a(u,x)\quad\text{such that}\quad \varphi(x,u)=f(x,a(u,x)),$$
where the operator function $F$ is defined by the Weyl (symmetric) symbol $f$ 
as a function of $m$ non-commuting variables $X_1,\ldots X_m$ [5:App.1].
\enddefinition

\remark{Remark 3} If the controls are absent and the pair $a(t)=(\frak A(t),
\bold e(t))$ is time-independent, Def.2 is reduced to the definition of the 
dynamical inverse problem of representation theory of the article [1].
\endremark

\remark{Remark 4} One may consider dynamical inverse problem of representation
theory for games, the interactively controlled systems and interactive games
(introduced by the author in [4]).
\endremark

\remark{Remark 5} One is able to interpret the correspondence 
$$\text{controlled system}\quad\longrightarrow\quad\text{representative 
dynamics}$$
as a quantization of the first. Such interpretation is very important for
the second quantization of intention fields in the interactive games [4].
\endremark

\remark{Remark 6} If the function $\varphi$ contained some constants
$c_\alpha\in\Bbb C$ then one may interpret them as time-independent variables
and include the matrices $C_\alpha\in\Mat_n(\Bbb C)$ instead of them in the
operator function $F$ (compare with the quantization of constants [2]).
\endremark

\Refs
\roster
\item"[1]" Juriev D.V., Dynamical inverse problem of representation theory
and non-commutative geometry [in Russian]. Fundam.Prikl.Matem. 4(1) (1998)
[e-version: funct-an/9507001].
\item"[2]" Juriev D.V., An excursus into the inverse problem of representation
theory [in Russian]. Report RCMPI-95/04 (1995) [e-version: mp\_arc/96-477].
\item"[3]" Vilenkin N.Ya., Special functions and group representation theory
[in Russian]. Moscow, Nauka, 1991.
\item"[4]" Juriev D., Interactive games and representation theory. I,II.
E-prints: math.FA/9803020, math.RT/9808098.
\item"[5]" Karasev M.V., Maslov V.P., Nonlinear Poisson brackets: geometry
and quantization. Amer.Math.Society, Providence, 1993.
\endroster
\endRefs
\enddocument